\documentclass[reqno,11pt]{amsart}
\usepackage{amsfonts}
\usepackage{latexsym,amssymb,amsthm,amsmath,graphics,amscd}

\theoremstyle{remark}

\theoremstyle{remark}

\theoremstyle{definition}

\numberwithin{equation}{section}

\begin{document}
\title[Yamabe solitons]{Conformal vector fields and Yamabe solitons}
\author[N. Turki]{Nasser Bin Turki}
\address{Department of Mathematics, College of science, King Saud University
P.O. Box-2455 Riyadh-11451, Saudi Arabia}
\email{nassert@ksu.edu.sa}
\author[B. Y. Chen]{Bang-Yen Chen}
\address{Department of Mathematics, Michigan State University, 619 Red Cedar
Road, East Lancing, Michigan 48824-1027 U.S.A.}
\email{chenb@msu.edu}
\author[S. Deshmukh]{Sharief Deshmukh}
\address{Department of Mathematics, College of science, King Saud University
P.O. Box-2455 Riyadh-11451, Saudi Arabia}
\email{shariefd@ksu.edu.sa}

\begin{abstract}
In this paper, we use less topological restrictions and more geometric and
analytic conditions to obtain some sufficient conditions on Yamabe solitons such that
their metrics are Yamabe metrics, that is, metrics of constant scalar curvature. 
More precisely, we use properties of conformal vector fields
to find several sufficient conditions on the soliton vector fields of Yamabe solitons
under which their metrics are of Yamabe metrics.
\end{abstract}

\maketitle

\section{\protect\medskip Introduction}

R. Hamilton introduced the notion of Yamabe flow (cf. \cite{9}), in which
the metric on a Riemannian manifold is deformed by evolving according to%
\begin{equation}\label{1.1}
\frac{\partial }{\partial t}g(t)=-R(t)g(t)\text{,}  
\end{equation}%
where $R(t)$ is the scalar curvature of the metric $g(t)$. Yamabe solitons
correspond to self-similar solutions of the Yamabe flow.

In dimension $n=2$ the Yamabe flow is equivalent to the Ricci flow (defined
by $\frac{\partial }{\partial t}g(t)=-2Ric(t)$, where $Ric$ denotes the
Ricci tensor). However, in dimension $n>2$ the Yamabe and Ricci flows do not
agree, since the first one preserves the conformal class of the metric but
the Ricci flow does not in general.

An $n$-dimensional connected Riemannian manifold $(M,g)$ is said to be a 
\textit{Yamabe soliton} if it admits a smooth vector field $\xi $ that
satisfies%
\begin{equation}\label{1.2}
\frac{1}{2}\pounds _{\xi }g=(R-\lambda )g\text{,}  
\end{equation}%
where $\pounds _{\xi }$ denotes the Lie derivative with respect to the
vector field $\xi $ and $\lambda $ is a constant. We call the vector field $%
\xi $ in the above definition a soliton vector field for $(M,g)$. In the
rest of this paper, we denote the Yamabe soliton satisfying \eqref{1.2} by $%
(M,g,\xi ,\lambda )$. A Yamabe soliton is said to be gradient Yamabe soliton
if the soliton vector field $\xi $ is gradient of a smooth function. Recall
that on a Riemannian manifold $(M,g)$, the metric $g$ is said to be a 
\textit{Yamabe metric} if the scalar curvature $R$ is a constant.

Yamabe flow and Yamabe solitons have been studied quite extensively (cf. 
\cite{1}-\cite{3}, \cite{5}, \cite{10}, \cite{11}). In \cite{6}, it is shown
that the metric of a compact gradient Yamabe soliton is a Yamabe metric and
the same result is achieved in \cite{10} giving a shorter proof. 

One of the interesting questions in the geometry of Yamabe solitons is to find
conditions on soliton vector fields of Yamabe solitons so that their metrics are of constant scalar curvature, that is, Yamabe metrics. There are two
options for such a study, such as imposing stronger topological restrictions
and less geometric and analytic conditions, or less topological restrictions
and more geometric and analytic conditions. In this paper, we study  
Yamabe solitons via the second option, merely with the topological
restriction that is connected. We then prove several new
sufficient conditions on the soliton fields of the Yamabe solitons for their metrics to be Yamabe metrics.

\section{\protect\medskip Preliminaries}

Recall that a smooth vector field $X$ on a Riemannian manifold $(M,g)$ is
said to be a conformal vector field if the local flow of $X$ consists of
local conformal transformations of $(M,g)$, which is equivalent to the fact
that the vector field $X$ satisfies%
\begin{equation*}
\frac{1}{2}\pounds _{\xi }g=\rho g\text{,}
\end{equation*}%
where $\rho $ is a smooth function on $M$ called the potential function of $%
X $. Thus, from the definition of Yamabe soliton $(M,g,\xi ,\lambda )$, it
follows that the soliton field $\xi $ is a conformal vector field with
potential function $R-\lambda $. Let $\eta $ be smooth $1$-form dual to the
soliton field $\xi $. Then, we define a skew symmetric tensor field $\varphi 
$ on the Yamabe soliton $(M,g,\xi ,\lambda )$ by%
\begin{equation}\label{2.1}
\frac{1}{2}d\eta (X,Y)=g(\varphi X,Y)\text{,\quad }X,Y\in \mathfrak{X}(M)%
\text{,}  \end{equation}%
where $\mathfrak{X}(M)$ is the Lie algebra of smooth vector fields on $M$.

Using Koszul's formula and equations \eqref{1.2}, \eqref{2.1}, the covariant derivative
of the soliton vector field $\xi $ is given by%
\begin{equation}\label{2.2}
\nabla _{X}\xi =(R-\lambda )X+\varphi X\text{,\quad }X\in \mathfrak{X}(M)%
\text{.}  \end{equation}%

Using equation \eqref{2.2}, a direct computation gives the following
expression for Riemannian curvature tensor of the Yamabe soliton $(M,g,\xi,\lambda )$:
\begin{equation}\label{2.3}
R(X,Y)\xi =X(R)Y-Y(R)X+\left( \nabla \varphi \right) (X,Y)-\left( \nabla
\varphi \right) (Y,X)\text{,}  
\end{equation}
where the covariant derivative $\left( \nabla \varphi \right)$ is defined by 
$$\left( \nabla \varphi \right) (X,Y)=\nabla_{X}\varphi Y-\varphi \left( \nabla _{X}Y\right).$$ 
Since the $2$-form $%
\Omega (X,Y)=g(\varphi X,Y)$ is closed, using equation \eqref{2.3}, a straight
forward computation leads to the following%
\begin{equation}\label{2.4}
\left( \nabla \varphi \right) (X,Y)=R(X,\xi )Y+Y(R)X-g(X,Y)\nabla R\text{%
,\quad }X,Y\in \mathfrak{X}(M)\text{,} \end{equation}%
where $\nabla R$ is the gradient of the scalar curvature $R.$ 

Let $Ric$ denote the Ricci tensor of the Yamabe soliton $(M,g,\xi ,\lambda )$. Then the Ricci
operator $Q$ of the Yamabe soliton is defined by $$g(QX,Y)=Ric(X,Y),$$ which
is a symmetric operator that satisfies%
\begin{equation}\label{2.5}
\sum \left( \nabla Q\right) (e_{i},e_{i})=\frac{1}{2}\nabla R\text{,} 
\end{equation}%
with $\left\{ e_{1},..,e_{n}\right\} $ being a local orthonormal frame on the
Yamabe soliton $(M,g,\xi ,\lambda )$. 

Using the skew-symmetry of the
operator $\varphi $, and equation \eqref{2.3}, we obtain%
\begin{equation*}
Ric(Y,\xi )=-(n-1)Y(R)-g\left( Y,\sum \left( \nabla \varphi \right)
(e_{i},e_{i})\right) \text{,}
\end{equation*}%
which leads to%
\begin{equation}
Q(\xi )=-(n-1)\nabla R-\sum \left( \nabla \varphi \right) (e_{i},e_{i})%
\text{.}  \label{2.6}
\end{equation}

Recall that the Laplace operator $\Delta $ acting on smooth vector fields on
a Riemannian manifold $(M,g)$ is an operator $\Delta :\mathfrak{X}%
(M)\rightarrow \mathfrak{X}(M)$ defined by%
\begin{equation*}
\Delta X=\sum \left( \nabla _{e_{i}}\nabla _{e_{i}}X-\nabla _{\nabla
_{e_{i}}}e_{i}X\right) \text{.}
\end{equation*}%
Thus, using equation \eqref{2.2}, we get the following%
\begin{equation}\label{2.7}
\Delta \xi =\nabla R+\sum \left( \nabla \varphi \right) (e_{i},e_{i})\text{.%
} \end{equation}

\section{Reducing metrics of Yamabe solitons to Yamabe metrics}

Now, we prove the first result of this paper.

\medskip

\textbf{Theorem 3.1. } \textit{Let $(M,g,\xi ,\lambda)$ be an $n$-dimensional connected Yamabe soliton with  $n>2$. If the soliton
field satisfies $\,Q(\xi )+\Delta \xi =\mu \xi\,$  for a positive
constant  $\mu $ and the scalar curvature $R$ of $M$ satisfies
\begin{equation*}
\lambda \leq R\leq \frac{n(n-1)}{n-2}\mu \text{,}
\end{equation*}%
then the metric $g$ is a Yamabe metric.}

\medskip

\begin{proof}
\medskip\ Using equations \eqref{2.6} and \eqref{2.7}, we have%
\begin{equation*}
Q(\xi )+\Delta \xi =-(n-2)\nabla R\text{,}
\end{equation*}%
which gives%
\begin{equation}\label{3.1}
\nabla R=-\frac{\mu }{(n-2)}\xi \text{.}  
\end{equation}%
Using above equation and equation \eqref{2.2}, we compute the Hessian operator $%
A_{R}$ of the scalar curvature $R$ to be%
\begin{equation*}
A_{R}(X)=-\frac{\mu (R-\lambda )}{n-2}X-\frac{\mu }{(n-2)}\varphi X\text{%
,\quad }X\in \mathfrak{X}(M)\text{,}
\end{equation*}%
that is,%
\begin{equation*}
A_{R}(X)+\frac{\mu (R-\lambda )}{n-2}X=-\frac{\mu }{(n-2)}\varphi X\text{.}
\end{equation*}%
The left-hand-side of above equation is symmetric, while the right-hand-side
is skew-symmetric and consequently we have%
\begin{equation}\label{3.2}
A_{R}(X)=-\frac{\mu (R-\lambda )}{n-2}X\text{,\quad }\varphi X=0\text{,\quad 
}X\in \mathfrak{X}(M)\text{.} 
\end{equation}%
Using the outcome $\varphi =0$, in the equation \eqref{2.6}, we obtain
$$Q(\xi)=-(n-1)\nabla R.$$ Taking divergence on both sides of this equation and
using equations \eqref{2.2}, \eqref{2.5}, with $\varphi =0$, we conclude%
\begin{equation*}
(R-\lambda )R+\frac{1}{2}\xi (R)=-(n-1)\Delta R\text{,}
\end{equation*}%
where $\Delta R=\rm{div}\left( \nabla R\right) $. Using equation \eqref{3.1},
the above equation leads to%
\begin{equation*}
(R-\lambda )R-\frac{n-2}{2\mu }\left\Vert \nabla R\right\Vert
^{2}=-(n-1)\Delta R\text{.}
\end{equation*}%
Now, using equation \eqref{3.2} in computing $\Delta R$, and inserting it in above
equation, we arrive at%
\begin{equation*}
\left( R-\lambda \right) \left( R-\frac{n(n-1)}{n-2}\mu \right) =\frac{n-2}{%
2\mu }\left\Vert \nabla R\right\Vert ^{2}\text{.}
\end{equation*}%
Finally, the bounds on scalar curvature together with above equation
confirms that the scalar curvature $R$ is a constant.
\end{proof}

\medskip

\textbf{Theorem 3.2. } \textit{Let }$(M,g,\xi ,\lambda )$\textit{\ be an }$n$%
-\textit{dimensional connected Yamabe soliton  with} $n>2$. \textit{If the soliton
field }$\xi $ \textit{annihilates the Hessian operator }$A_{R}$ \textit{and} 
\textit{satisfies }$\Delta \xi =-\mu \xi $\textit{ for a positive constant }%
$\mu $ \textit{and if the Ricci curvature in the direction of the soliton field
is a constant,\ then the metric }$g$ \textit{is a Yamabe metric.}

\begin{proof}
As the Ricci curvature in the direction of $\xi $ is a constant, we have $%
Ric(\xi ,\xi )=\nu \left\Vert \xi \right\Vert ^{2}$, where $\nu $ is a
constant. Using equations \eqref{2.6}, \eqref{2.7} and $\Delta \xi =-\mu \xi $, we get%
\begin{equation*}
Q(\xi )=\mu \xi -(n-2)\nabla R\text{,}
\end{equation*}%
which on taking the inner product with $\xi $ leads to%
\begin{equation*}
\left( \mu -\nu \right) \left\Vert \xi \right\Vert ^{2}=(n-2)g(\xi ,\nabla R)%
\text{.}
\end{equation*}%
Taking covariant derivative with respect to $X$ in above equation and using
equation \eqref{2.2} and the fact that $\xi $ annihilates $A_{R}$, we conclude%
\begin{equation*}
2\left( \mu -\nu \right) g\left( (R-\lambda )X+\varphi X,\xi \right)
=(n-2)g\left( (R-\lambda )X+\varphi X,\nabla R\right) \text{,}
\end{equation*}%
that is,%
\begin{equation*}
2\left( \mu -\nu \right) (R-\lambda )\xi -(n-2)\nabla R=\varphi \left(
2\left( \mu -\nu \right) (R-\lambda )\xi -(n-2)\nabla R\right) \text{.}
\end{equation*}%
Taking the inner product with $2\left( \mu -\nu \right) (R-\lambda )\xi
-(n-2)\nabla R$ in above equation and using the skew-symmetry of $\varphi $,
we conclude that%
\begin{equation}
\nabla R=2\frac{\left( \mu -\nu \right) (R-\lambda )}{n-2}\xi =h\xi \text{,}
\label{3.3}
\end{equation}%
where%
\begin{equation}
h=2\frac{\left( \mu -\nu \right) (R-\lambda )}{n-2}\text{.}  \label{3.4}
\end{equation}%
Taking covariant derivative with respect to $X$ in the equation \eqref{3.3} and
using equations \eqref{2.2}, we get%
\begin{equation*}
A_{R}X=2\frac{\left( \mu -\nu \right) }{n-2}X(R)\xi +h(R-\lambda )X+h\varphi
X\text{,}
\end{equation*}%
which on using the symmetry of operator $A_{R}$ and the skew-symmetry of $%
\varphi $ and equation \eqref{3.3}, leads to%
\begin{eqnarray*}
2h\varphi X &=&2\frac{\left( \mu -\nu \right) }{n-2}\left( g(X,\xi )\nabla
R-g(X,\nabla R)\xi \right) \\
&=&4\left( \frac{\mu -\nu }{n-2}\right) ^{2}(R-\lambda )\left( g(X,\xi )\xi
-g(X,\xi )\xi \right) =0\text{.}
\end{eqnarray*}%
Thus, we have either $h=0$ or $\varphi =0$. If $h=0$, then equation \eqref{3.3}
implies $R$ is a constant. Therefore, suppose $\varphi =0$, which in view of
equation \eqref{2.6}, gives%
\begin{equation*}
Q(\xi )=-(n-1)\nabla R\text{.}
\end{equation*}%
Taking the inner product with $\xi $ in above equation implies%
\begin{equation*}
\nu \left\Vert \xi \right\Vert ^{2}=-(n-1)\xi (R)\text{,}
\end{equation*}%
which on using equation \eqref{3.3}, gives%
\begin{equation*}
\nu \left\Vert \xi \right\Vert ^{2}=-(n-1)h\left\Vert \xi \right\Vert ^{2}%
\text{,}
\end{equation*}%
that is,%
\begin{equation*}
\left\Vert \xi \right\Vert ^{2}\left( 2\frac{\left( \mu -\nu \right)
(R-\lambda )(n-1)}{n-2}+\nu \right) =0\text{.}
\end{equation*}%
Note that if $\mu =\nu $, then equation \eqref{3.3} would imply $R$ is a constant.
Also, if $\xi =0$, the equation \eqref{1.2} will imply $R=\lambda $, therefore,
now suppose $\mu \neq \nu $, then the above equation gives%
\begin{equation*}
R=\lambda -\frac{(n-2)\nu }{2(n-1)\left( \mu -\nu \right) }\text{,}
\end{equation*}%
that is, $R$ is a constant.
\end{proof}

\medskip 

In the above theorem, we could replace the condition that soliton field $\xi 
$ annihilates the Hessian operator $A_{R}$ in exchange with condition on
Ricci curvature in the direction of $\xi $ being constant is made stronger,
namely $\xi $ satisfies $Q(\xi )=\nu \xi $ for a constant $\nu $ and certain
bounds on the scalar curvature. We get the following result.

\medskip

\textbf{Proposition 3.1. } \textit{Let }$(M,g,\xi ,\lambda )$\textit{\ be an 
}$n$-\textit{dimensional connected Yamabe soliton }$n>2$. \textit{If the
soliton field }$\xi $ \textit{satisfies }$\Delta \xi =-\mu \xi $\textit{,
and }$Q(\xi )=\nu \xi $ \textit{for constants }$\mu \geq 0$ \textit{and }$%
\nu $, and scalar curvature satisfies%
\begin{equation*}
\lambda \leq R\leq n(n-1)\mu \text{,}
\end{equation*}%
\textit{then the metric }$g$ \textit{is a Yamabe metric.}

\begin{proof}
Using equations \eqref{2.6} and \eqref{2.7}, we have $$Q(\xi )+\Delta \xi =-(n-2)\nabla R$$
 which gives
\begin{equation}
\nabla R=\frac{1}{n-2}(\mu -\nu )\xi \text{.}  \label{3.5}
\end{equation}%
Taking covariant derivative in above equation with respect to $X$ and using
equation \eqref{2.2} gives%
\begin{equation*}
A_{R}X-\frac{1}{n-2}(\mu -\nu )(R-\lambda )X=\frac{1}{n-2}(\mu -\nu )\varphi
X\text{,}
\end{equation*}%
which has one side symmetric while other skew-symmetric, and therefore
yields $\varphi =0$. Then equations \eqref{2.6} and \eqref{2.7} respectively give $$\nu
\xi =-(n-1)\nabla R \quad {\rm and}\quad \mu \xi =-\nabla R,$$ which lead to%
\begin{equation}
\nu =(n-1)\mu \text{.}  \label{3.6}
\end{equation}%
Thus, inserting the value of $\nu $ in equation \eqref{3.5} and in $Q(\xi )=\nu
\xi $ respectively, give%
\begin{equation}
\nabla R=-\mu \xi \text{,\quad }Q(\xi )=(n-1)\mu \xi \text{.}  \label{3.7}
\end{equation}%
Note that if $\mu =0$, then by above equation $R$ is a constant. Therefore from
now on we will assume $\mu >0$. Taking divergence in the second equation of
\eqref{3.7} and using equations \eqref{2.2} and \eqref{2.5} we get%
\begin{equation*}
R(R-\lambda )+\frac{1}{2}\xi (R)=n(n-1)\mu (R-\lambda )
\end{equation*}%
and multiplying above equation by $-\mu $ and using first equation in \eqref{3.7}
yields%
\begin{equation*}
\mu (R-\lambda )\left( n(n-1)\mu -R\right) +\frac{1}{2}\left\Vert \nabla
R\right\Vert ^{2}=0\text{.}
\end{equation*}%
Hence, using bounds on scalar curvature in above equation implies that $R$
is a constant.
\end{proof}

Using the motivation of analytic vector fields on complex manifolds whose local flows
leave complex structure invariant,  the notion of $\varphi $%
-analytic conformal vector fields was introduced in \cite{7}. According to \cite{7}, a {\it $\varphi $-analytic
conformal vector field} is a conformal vector field with respect whose local
flow, the tensor field $\varphi $ appearing in equation \eqref{2.2} is invariant.
For examples of $\varphi $-analytic conformal vector fields, see \cite{8}.

Equations \eqref{3.1} and \eqref{3.3} in the proofs of Theorems 3.1 and 3.2, and the first
equation in \eqref{3.7} suggest that soliton fields under given conditions are
 $\varphi $-analytic conformal vector fields (cf. \cite{7}).  In
the rest of this paper, we will use this property that the soliton fields being $\varphi $-analytic to
prove following two results. First, observe that if the soliton field of the
Yamabe soliton $(M,g,\xi ,\lambda )$ is $\varphi $-analytic, then we have%
\begin{equation*}
\pounds _{\xi }\varphi =0\text{,}
\end{equation*}%
which in view of equation \eqref{2.2} is equivalent to%
\begin{equation}
\left( \nabla \varphi \right) (\xi ,X)=0\text{,\quad }X\in \mathfrak{X}(M)%
\text{.}  \label{3.8}
\end{equation}%
It is known that the soliton field $\xi $ is $\varphi $-analytic if and
only if $\nabla R=f\xi $ for a smooth function $f$ on $M$ (cf. \cite{7}).
Moreover, if $f=c$, a constant, then $\xi $ is called a {\it constant type-$c$} $%
\varphi $-analytic conformal vector field.

\medskip

\textbf{Theorem 3.3. } \textit{Let }$(M,g,\xi ,\lambda )$\textit{\ be an }$n$%
-\textit{dimensional connected Yamabe soliton.} \textit{If the soliton field
is a }$\varphi $-\textit{analytic vector field of constant type-}$c$ {\it with} $c<0$ 
\textit{and the scalar curvature satisfies}%
\begin{equation*}
\lambda \leq R\leq -n(n-1)c\text{,}
\end{equation*}%
\textit{\ then the metric }$g$ \textit{is a Yamabe metric.}

\begin{proof}
Since $\xi $ is $\varphi $-analytic vector field of constant type $c$, we
have%
\begin{equation*}
\nabla R=c\xi \text{,}
\end{equation*}%
which in view of equation \eqref{2.2}, gives%
\begin{equation*}
A_{R}X=c(R-\lambda )X+c\varphi X\text{.}
\end{equation*}%
Using the symmetry of $A_{R}$ and the skew-symmetry of $\varphi $, we get $%
\varphi =0$. Now, using equation \eqref{2.4}, we get%
\begin{equation*}
R(X,\xi )Y=g(X,Y)\nabla R-Y(R)X\text{,}
\end{equation*}%
which in view of $\nabla R=c\xi $, gives%
\begin{equation*}
Ric(Y,\xi )=-(n-1)cg(Y,\xi ).
\end{equation*}%
Thus,%
\begin{equation}
Q(\xi )=-(n-1)c\xi \text{.}  \label{3.9}
\end{equation}

Using equations \eqref{2.5} and \eqref{2.2}, a straight forward computation gives 
$$\rm{div}\left( Q(\xi )\right) =R(R-\lambda )+\frac{1}{2}\xi (R)\quad {\rm and}\quad \rm{div}%
\,\xi =n(R-\lambda ).$$ Now, taking divergence in equation \eqref{3.9} and
multiplying by $c$ gives%
\begin{equation*}
cR(R-\lambda )+\frac{1}{2}\left\Vert \nabla R\right\Vert
^{2}=-n(n-1)c^{2}(R-\lambda )\text{,}
\end{equation*}%
that is,%
\begin{equation*}
c(R-\lambda )\left( R+n(n-1)c\right) +\frac{1}{2}\left\Vert \nabla
R\right\Vert ^{2}=0\text{.}
\end{equation*}%
Then the bounds on the scalar curvature prove that $R$ is a constant.
\end{proof}

\medskip

The following theorem uses the combination of the length of soliton vector
field (the speed of the integral curves of the soliton field) and scalar
curvature, that is, the speed of the integral curves of soliton field is
controlled by the scalar curvature.

\medskip

\noindent \textbf{Theorem 3.4.} \textit{Let }$(M,g,\xi ,\lambda )$\textit{\
be a connected Yamabe soliton of positive scalar curvature. If the soliton
field is }$\varphi $\textit{-analytic vector field and the length of soliton
field satisfies}%
\begin{equation*}
\left\Vert \xi \right\Vert =\sqrt{2R}\text{,}
\end{equation*}%
\textit{\ then }$g$\textit{\ is a Yamabe metric.}

\begin{proof}
Using $\left\Vert \xi \right\Vert ^{2}=2R$ and equation \eqref{2.2}, we get%
\begin{equation}
\nabla R=(R-\lambda )\xi -\varphi \xi \text{,}  \label{3.10}
\end{equation}%
which on taking the inner product with $\xi $, leads to%
\begin{equation}
\xi (R)=2(R-\lambda )R\text{.}  \label{3.11}
\end{equation}%
Now, multiplying equation \eqref{3.10} by $2R$ gives%
\begin{equation}
2R(R-\lambda )\xi -2R\varphi \xi =\left\Vert \xi \right\Vert ^{2}\nabla R%
\text{.}  \label{3.12}
\end{equation}%
Since $\xi $ is $\varphi $-analytic, after using equation \eqref{3.5} and equation
\eqref{2.4}, we find $X(R)\xi =g(X,\xi )\nabla R$, which on substitution $X=\xi $
reads%
\begin{equation}
\xi (R)\xi =\left\Vert \xi \right\Vert ^{2}\nabla R\text{.}  \label{3.13}
\end{equation}%
Inserting this equation in equation \eqref{3.12} gives $$2R(R-\lambda )\xi
-2R\varphi \xi =\xi (R)\xi ,$$ which in view of equation \eqref{3.11} and $R>0$
yields
\begin{equation}
\varphi \xi =0\text{.}  \label{3.14}
\end{equation}%
Now, $\xi $ being a $\varphi $-analytic vector field, we have $\nabla R=f\xi 
$ for a smooth function $f$, which gives%
\begin{equation*}
A_{R}X=X(f)\xi +f(R-\lambda )X+f\varphi X\text{.}
\end{equation*}%
Using the symmetry of $A_{f}$ and the skew-symmetry of $\varphi $ in above
equation, we obtain
\begin{equation}
2f\varphi X=g(X,\xi )\nabla f-X(f)\xi \text{,}  \label{3.15}
\end{equation}%
which on using $X=\xi $ and equation \eqref{3.14} yields $\left\Vert \xi
\right\Vert ^{2}\nabla f=\xi (f)\xi $, that is,%
\begin{equation}
\left\Vert \xi \right\Vert ^{2}\left\Vert \nabla f\right\Vert ^{2}=\xi
(f)^{2}\text{.}  \label{3.16}
\end{equation}%
Now, taking a local orthonormal frame $\left\{ e_{1},e_{2},..e_{n}\right\} $%
, the equation \eqref{3.15} gives%
\begin{equation*}
4f^{2}\left\Vert \varphi \right\Vert ^{2}=4f^{2}\sum \left\Vert \varphi
e_{i}\right\Vert ^{2}=2\left\Vert \xi \right\Vert ^{2}\left\Vert \nabla
f\right\Vert ^{2}-2\xi (f)^{2}\text{,}
\end{equation*}%
which in view of equation \eqref{3.16}, gives%
\begin{equation*}
f^{2}\left\Vert \varphi \right\Vert ^{2}=0\text{.}
\end{equation*}%
If $f=0$, then $R$ is a constant. Therefore, suppose $\varphi =0$. Then the
equation \eqref{3.10} takes the form $\nabla R=(R-\lambda )\xi $, taking
divergence on both sides of this equation, we get%
\begin{equation}
\Delta R=\xi (R)+n(R-\lambda )^{2}\text{.}  \label{3.17}
\end{equation}%
Now, using $\varphi =0$ in equation \eqref{2.2}, it takes the form $\nabla _{X}\xi
=(R-\lambda )X$ and using it for finding the divergence of $Q(\xi )$, gives $%
\rm{div}\,Q(\xi )=R(R-\lambda )+\frac{1}{2}\xi (R)$, which in view of
equation \eqref{3.11} implies%
\begin{equation}
\rm{div}\, Q(\xi )=2R(R-\lambda )\text{.}  \label{3.18}
\end{equation}%
Also, equation \eqref{2.6} gives $Q(\xi )=-(n-1)\nabla R$, which on taking
divergence both sides and using equation \eqref{3.18} implies%
\begin{equation*}
\Delta R=-\frac{2}{n-1}R(R-\lambda )\text{.}
\end{equation*}%
Above equation together with equations \eqref{3.11} and \eqref{3.17} gives%
\begin{equation*}
n(R-\lambda )\left( \frac{2}{n-1}+R-\lambda \right) =0
\end{equation*}%
and it proves $R$ is a constant
\end{proof}

\medskip

\section{Further results}

In this section, we analyze the impact on the Yamabe
soliton if the soliton field $\xi $ is closed as well as if it is not
closed. Recall that by equation \eqref{2.2}, if the soliton field is closed, then $%
\varphi =0$ and if it is not closed, then $\varphi \neq 0$. We shall call
the operator $\varphi $ the operator associated to soliton field.

\medskip

\noindent \textbf{Theorem 4.1.} \textit{Let }$(M,g,\xi ,\lambda )$\textit{\
be an }$n$-\textit{dimensional  connected Yamabe soliton. If the soliton
field }$\xi $\textit{\ is closed and the scalar curvature }$R$\textit{\ is a
constant along the integral curves of }$\xi $\textit{, then the metric }$g$ 
\textit{is a Yamabe metric.}

\begin{proof}
Let $\xi $ be closed. Then equation \eqref{2.2} takes the form%
\begin{equation*}
\nabla _{X}\xi =\left( R-\lambda \right) X\text{.}
\end{equation*}%
Define $f=\frac{1}{2}\left\Vert \xi \right\Vert ^{2}$, which on using above
equation implies the following expression for the gradient $\nabla
f=(R-\lambda )\xi $. Thus the Hessian operator $A_{f}$ of the function $f$
is given by%
\begin{equation*}
A_{f}(X)=X(R)\xi +\left( R-\lambda \right) ^{2}X\text{.}
\end{equation*}%
Now, taking the inner product in above equation with a vector field $Y$ and
using the fact that the Hessian operator is symmetric we arrive at%
\begin{equation*}
X(R)g(Y,\xi )=Y(R)g(X,\xi )\text{,}
\end{equation*}%
which leads to%
\begin{equation*}
X(R)\xi =g(X,\xi )\nabla R\text{.}
\end{equation*}%
Taking the inner product in above equation with $\nabla R$ and replacing $X$
by $\xi $ leads to%
\begin{equation*}
\left\Vert \xi \right\Vert ^{2}\left\Vert \nabla R\right\Vert ^{2}=\xi
(R)^{2}\text{.}
\end{equation*}%
As the scalar curvature $R$ is a constant along the integral curves of $\xi $%
, that is, $\xi (R)=0$, above equation implies that $R$ is a constant.
\end{proof}

\medskip 
Now, in contrast to above theorem, we prove the following
in case that the soliton field $\xi $ is not closed.\medskip 

\noindent \textbf{Theorem 4.2.} \textit{Let }$(M,g,\xi ,\lambda )$\textit{\
be an }$n$-\textit{dimensional connected Yamabe soliton. If the soliton
field }$\xi $\textit{\ is not closed and it annihilates the associated
operator }$\varphi $\textit{, then the metric }$g$ \textit{is a Yamabe
metric.}

\begin{proof}
Suppose that the soliton field $\xi $ is not closed (that is, $\varphi \neq 0
$) and that it annihilates $\varphi $, that is, $\varphi (\xi )=0$. Then the
function $f=\frac{1}{2}\left\Vert \xi \right\Vert ^{2}$ on using equation
\eqref{2.2}, has for its gradient%
\begin{equation*}
\nabla f=\left( R-\lambda \right) \xi \text{.}
\end{equation*}%
Using equation \eqref{2.2} and above equation, we find the following expression
for the Hessian operator $A_{f}$%
\begin{equation*}
A_{f}(X)=X(R)\xi +\left( R-\lambda \right) ^{2}X+(R-\lambda )\varphi X\text{.%
}
\end{equation*}%
Now, as the Hessian operator $A_{f}$ is symmetric and the associated
operator $\varphi $ symmetric, the above equation gives%
\begin{equation}
(R-\lambda )\varphi X=\frac{1}{2}\left( g(X,\xi )\nabla R-X(R)\xi \right) 
\text{,}  \label{4.1}
\end{equation}%
which on taking $X=\xi $ leads to%
\begin{equation*}
\xi (R)\xi =\left\Vert \xi \right\Vert ^{2}\nabla R\text{.}
\end{equation*}%
Taking the inner product in above equation with $\nabla R$ gives%
\begin{equation}
\xi (R)^{2}=\left\Vert \xi \right\Vert ^{2}\left\Vert \nabla R\right\Vert
^{2}\text{.}  \label{4.2}
\end{equation}%
Also, using equation \eqref{4.1} and a local orthonormal frame $\left\{
e_{1},..,e_{n}\right\} $, we compute%
\begin{equation*}
(R-\lambda )^{2}\left\Vert \varphi \right\Vert
^{2}=\sum\limits_{i}\left\Vert (R-\lambda )\varphi e_{i}\right\Vert ^{2}=%
\frac{1}{4}\left( 2\left\Vert \xi \right\Vert ^{2}\left\Vert \nabla
R\right\Vert ^{2}-2\xi (R)^{2}\right) \text{,}
\end{equation*}%
which in view of equation \eqref{4.2} leads to%
\begin{equation*}
(R-\lambda )^{2}\left\Vert \varphi \right\Vert ^{2}=0\text{.}
\end{equation*}%
Since $\xi $ is not closed, that is, $\varphi \neq 0$, we find that $R$ is a
constant.
\end{proof}

\medskip 

\textbf{Acknowledgements:} This work is supported by King Saud University,
Deanship of Scientific Research, College of Science Research Center.

\end{document}